\documentclass[12pt]{amsart}
\usepackage{amssymb}
\usepackage{amsfonts}
\usepackage{amscd}
\usepackage[all]{xypic}

\swapnumbers

\def\ord{{\rm ord}}

\DeclareMathOperator*{\sq}{\square}
\DeclareMathOperator*{\wedgem}{\wedge}

\let\cal\mathcal

\def\11{{\mathbf 1}}
\def\AA{{\mathbf A}}

\def\CC{{\mathbf C}}

\def\FF{{\mathbf F}}
\def\GG{{\mathbf G}}

\def\LL{{\mathbf L}}

\def\NN{{\mathbf N}}

\def\QQ{{\mathbf Q}}

\def\ZZ{{\mathbf Z}}

\def\cA{{\mathcal A}}
\def\cB{{\mathcal B}}
\def\cC{{\mathcal C}}

\def\cL{{\mathcal L}}

\def\cO{{\mathcal O}}

\def\cT{{\mathcal T}}

\mathchardef\alphag="7C0B \mathchardef\betag="7C0C
\mathchardef\gammag="7C0D \mathchardef\deltag="7C0E
\mathchardef\varepsilong="7C22 \mathchardef\varphig="7C27
\mathchardef\psig="7C20 \mathchardef\zetag="7C10
\mathchardef\epsilong="7C0F \mathchardef\rhog="7C1A
\mathchardef\taug="7C1C \mathchardef\upsilong="7C1D
\mathchardef\iotag="7C13 \mathchardef\thetag="7C12
\mathchardef\pig="7C19 \mathchardef\sigmag="7C1B
\mathchardef\etag="7C11 \mathchardef\omegag="7C21
\mathchardef\kappag="7C14 \mathchardef\lambdag="7C15
\mathchardef\mug="7C16 \mathchardef\xig="7C18
\mathchardef\chig="7C1F \mathchardef\nug="7C17
\mathchardef\varthetag="7C23 \mathchardef\varpig="7C24
\mathchardef\varrhog="7C25 \mathchardef\varsigmag="7C26
\mathchardef\Omegag="7C0A \mathchardef\Thetag="7C02
\mathchardef\Sigmag="7C06 \mathchardef\Deltag="7C01
\mathchardef\Phig="7C08 \mathchardef\Gammag="7C00
\mathchardef\Psig="7C09 \mathchardef\Lambdag="7C03
\mathchardef\Xig="7C04 \mathchardef\Pig="7C05
\mathchardef\Upsilong="7C07

\newtheorem{theorem}[subsection]{Theorem}
\newtheorem{lem}[subsection]{Lemma}

\newtheorem{prop}[subsection]{Proposition}

\newtheorem{conjecture}[subsection]{Conjecture}

\theoremstyle{definition}
\newtheorem{definition}[subsection]{Definition}

\newtheorem{def-prop}[subsection]{Proposition-Definition}
\newtheorem{def-theorem}[subsection]{Theorem-Definition}
\newtheorem{def-lem}[subsection]{Lemma-Definition}

\theoremstyle{remark}
\newtheorem{remark}[subsection]{Remark}

\theoremstyle{plain}

\numberwithin{equation}{section}

\def\boxit#1#2{\setbox1=\hbox{\kern#1{#2}\kern#1}%
\dimen1=\ht1 \advance\dimen1 by #1 \dimen2=\dp1 \advance\dimen2 by
#1
\setbox1=\hbox{\vrule height\dimen1 depth\dimen2\box1\vrule}%
\setbox1=\vbox{\hrule\box1\hrule}%
\advance\dimen1 by .4pt \ht1=\dimen1 \advance\dimen2 by .4pt
\dp1=\dimen2 \box1\relax}

\renewcommand{\theequation}{\thesubsection.\arabic{equation}}
\def\ord{{\rm ord}}
\def\ac{{\overline{\rm ac}}}
\def\LPas{\cL_{\rm DP}}
\def\Ltame{\cL_{\rm DP}^{\rm tame}}
\def\To{\cT_0}
\def\Tinf{\cT_{\infty}}
\def\Tinfd{\cT_{\infty}^{(d)}}

\def\Gal{{\rm Gal}}

\def\llp{\mathopen{(\!(}}

\def\rrp{\mathopen{)\!)}}

\def\chara{{\rm char}}
\DeclareMathOperator*{\Spec}{Spec}

\author{Raf Cluckers}
\address{Katholieke Universiteit Leuven, Departement wiskunde,
Celestijnenlaan 200B, B-3001 Leu\-ven, Bel\-gium. Current address:
\'Ecole Normale Sup\'erieure, D\'epartement de
ma\-th\'e\-ma\-ti\-ques et applications, 45 rue d'Ulm, 75230 Paris
Cedex 05, France} \email{cluckers@ens.fr}
\urladdr{www.dma.ens.fr/$\sim$cluckers/}

\author{Jan Denef}
\address{Katholieke Universiteit Leuven, Departement wiskunde,
Celestijnenlaan 200B, B-3001 Leu\-ven, Bel\-gium.}
\email{jan.denef@wis.kuleuven.be}
\urladdr{www.wis.kuleuven.be/algebra/denef.html}



\title[]{Orbital integrals for linear groups}

\begin{document}

\begin{abstract}
For a linear group $G$ acting on an absolutely irreducible variety
$X$ over $\QQ$, we describe the orbits of $X(\QQ_p)$ under
$G(\QQ_p)$ and of $X(\FF_p\llp t\rrp)$ under $G(\FF_p\llp t\rrp)$
for $p$ big enough. This allows us to show that the degree of a
wide class of orbital integrals over $\QQ_p$ or $\FF_p\llp t\rrp$
is $\leq 0$ for $p$ big enough, and similarly for all finite field
extensions of $\QQ_p$ and $\FF_p\llp t\rrp$.
\end{abstract}

\maketitle

\begin{enumerate}
\item[1.] Introduction  \item[2.] The logical setting and pseudo-finite
fields \item[3.] Cohomological lemmas \item[4.] Definability of
cohomology
\item[5.] Proof of main results
\end{enumerate}

\section{Introduction}

Let $F$ be a number field with ring of integers $\cO_F$. Let $\cA_F$
be the collection of all finite field extensions of non-archimedean
  completions of $F$. Let $\cB_F$ be the collection of
all fields of the form $\FF_q\llp t\rrp$ which are rings over
$\cO_F$. For $K$ in $\cA_F\cup\cB_F$, let $\cO_K$ be its valuation
ring, $M_K$ its maximal ideal,
 $k_K$
  its residue field, and $q_K:=\sharp k_K$. For $N>0$,
let $\cC_N$ be
$$
\cC_N:=\{K\in \cA_F\cup\cB_F\mid \chara(k_K)>N\}.
$$
For any $K\in\cC_1$, let $\ac:K^\times\to k_ K^\times$ be a
multiplicative map extending the projection $\cO_K^\times\to
k_K^\times$, put $\ac(0)=0$, and let $\ord: K^\times\to\ZZ$ be the
order.

Let $G$ be a linear algebraic group over $F$, rationally acting on
an absolutely irreducible\footnote{See Remark \ref{remirr} to loosen
the condition of absolute irreducibility of $X$.} algebraic variety
$X$ over $F$. Suppose that $X$ is a homogeneous $G$-space, that is,
the action of $G(\CC)$ on $X(\CC)$ is transitive. For $K$ a field
over $F$ and $x\in X(K)$, let $G(K)(x)$ be the orbit of $x$ under
the action of $G(K)$.

\begin{theorem}\label{mt1}Let $U\subset X$ be an affine open.
Then there are finitely many regular functions $f_i:U\to \AA^1_F$
and integers $N>0$ and $d>0$ such that, for any $K\in \cC_N$ and
any $x\in X(K)$, the set $G(K)(x)\cap U(K)$ depends only on
$\ac(f_i(x))$ and $\ord(f_i(x))\bmod d$.
\end{theorem}

\begin{definition}\label{IKx}
Let $U\subset X$ be an affine open and let $f,g_i:U\to \AA^1_F$ be
regular functions. Let $\omega$ be a volume form on $U$, that is, a
degree $n$ rational differential form on $U$ when $X$ is of
dimension $n$. For each $K\in \cC_1$ let $W(K)$ be $\cap_i
g^{-1}_i(U_i)$ with $U_i$ either $\cO_K$ or $M_K$. For each $K\in
\cC_1$ and for each $x$ in $X(K)$, under the condition of
integrability for all $s>0$, consider the orbital integral
\begin{equation}\label{or1}
I_{K,x}(s):=\int_{G(K)(x) \cap W(K)} |f|^s\, |\omega|_K,
\end{equation}
with $|\omega|_K$ the measure on $U(K)$ associated to $\omega$. If
for some $K$ and $x$ this is not integrable, put $ I_{K,x}(s):=0$
for this $K$ and $x$.
\end{definition}

\begin{theorem}\label{mt3}
Let $I_{K,x}(s)$ be as in Definition \ref{IKx}. Then, there exists
$N>0$ such that $I_{K,x}(s)$ is rational and of degree $\leq 0$ in
$q_K^{-s}$ for each $K\in\cC_N$ and for each $x\in X(K)$.
\end{theorem}
 \bigskip
\noindent {\bf Addendum to Theorem \ref{mt3}.}\textit{ Let
$I_{K,x}(s)$ be as in Definition \ref{IKx}. Then there exist $N>0$
and $F(\LL,T)$ in $\QQ(\LL,T)$ of the form
$$
F(\LL,T)=T^{a}\, \prod_{i=1}^m (1-\LL^{a_i}T^{b_i}),
$$
with $a,a_i,b_i\in\ZZ$, $b_i> 0$, and $m\geq 0$, such that
$$
F(q_K,q_K^{-s})I_{K,x}(s)
$$
is a polynomial in $q_K^{-s}$ for each $K$ in $\cC_N$ and each $x\in
X(K)$. In particular, by Theorem \ref{mt3}, the degrees of the
numerator and the denominator of $I_{K,x}(s)$ are uniformly bounded
and only finitely many real poles with bounded multiplicities in $s$
can occur in $I_{K,x}(s)$ when $K$ varies over $\cC_N$ for suitable
$N$ and when $x$ varies in $X(K)$.}

\bigskip
\noindent
\subsection{}
Orbital integrals as in Theorem \ref{mt3} occur in representation
theory and in the study by Igusa \cite{Igusa:intro}, \cite{Igusa4}
and others, of prehomogeneous vector spaces, where $X=\AA_F^n$, $G$
acts linearly on $X$, $f$ is a relative invariant of this action,
and $W(K)$ is a Cartesian product of sets of the form $\cO_K$ and
$a+M_K$ with $a\in \cO_K$.\footnote{Note that the characteristic
function of such $W(K)$ is an important kind of Schwartz-Bruhat
function, since it is the essential building block for
Schwartz-Bruhat functions on adeles and ideles.}
 For such $X$, $f$,
and $W(K)$, Igusa \cite{Igusa4} determines, under some extra
conditions, the poles of the integral (\ref{or1}) in terms of
explicit group theoretical invariants. Since the multiplicity of the
poles is a priori bounded by the dimension of $X$, this gives an
explicit bound on the degree of the denominator of (\ref{or1}).
Combining this with Theorem \ref{mt3} then also gives a bound on the
degree of the numerator in terms of these group theoretical data.

For general $f$, Theorem \ref{mt3} and its addendum give a single
finite set of candidate poles in $t=q_K^{-s}$ of (\ref{or1}) and
uniform bounds on the degree of the denominator and the numerator of
(\ref{or1}) when $K$ varies in $\cC_N$, generalizing greatly the
situation of \cite{Igusa4}.

\section{The logical setting and pseudo-finite fields}\label{sec:setting}

\subsection{The languages $\LPas$ and $\Ltame$}\label{slt}

Let $\LL_{{\rm Ord}}=(+,-,\leq,0)$ be the language of ordered
groups.

For $K$ a valued field, $M$ the maximal ideal of its valuation
ring $R$, and $k$ the residue field, an \emph{angular component
modulo $M$} (or angular component for short) is a map $\ac:K\to k$
such that the restriction to $K^\times$ is a multiplicative
homomorphism to $k^\times$, the restriction to $R^\times$
coincides with the restriction to $R^\times$ of the natural
projection $R\to k$, and such that $\ac(0)=0$.

The \emph{language $\LPas$ of Denef-Pas} is defined as the three
sorted language
 $$(\LL_{\rm Rings},\LL_{{\rm Rings}},\LL_{\rm Ord},
 \ac,\ord).$$

The sorts of $\LPas$ are a valued field $K$ with residue field $k$
and value group $G$. The function symbol $\ord$ is the additively
written valuation $\ord:K^\times\to G$, and $\ac:K\to k$ is an
angular component.\footnote{The problem that the function $\ord$ is
not defined globally on $K$ is easily settled and the reader may
choose a way to do so. For example, the reader may choose a value of
$\ord(0)$ in the value group and treat the cases that the argument
of $\ord$ equals zero always separately, or, the reader may add a
symbol $+\infty$ to the language $\LL_{\rm Ord}$ that is bigger than
any element of the value group, and make the natural changes. We
will always make clear what we mean by expressions as $\ord(x)\leq
\ord(y)$ and so on when $x$ or $y$ can be zero.}
 The first ring
language is used for the valued field, the second for the residue
field, and $\LL_{\rm Ord}$ is used for the value group.

The language $\Ltame$ is obtained from $\LPas$ by removing the value
group sort and replacing it by infinitely many sorts for the
quotients $G/nG$, $n=2,3,\ldots$ where $G$ is the value group. Here,
$G/nG$ is considered as a group and the language for each of these
sorts is the group language $(+,-,0)$ together with projection maps
$\pi_{nm}:G/nG\to G/mG$ for $m\geq 2$ a divisor of $n$ and order
maps $\ord_n:K^\times\to G/nG$ making commuting diagrams, extended
by $\ord_n(0)=0$.

\subsection{Theories $\To$, $\Tinf$, and $\Tinfd$}\label{sstt}
In all what follows $\cT$ is any theory in $\LPas$ that contains
the $\LPas$-theory $\To$ of Henselian valued fields with angular
component modulo the maximal ideal and with residual
characteristic zero. This theory has elimination of valued field
quantifiers in $\LPas$ by \cite{Pas}.
 Recall that a
  perfect
   field $k$ is called pseudo-finite if it has a unique
field extension of any given finite degree and if any absolutely
irreducible variety over $k$ has a $k$-rational point.
 The theory of pseudo-finite fields is a first order theory which can
be expressed by an infinite axiom scheme in the language of rings.

Let $\Tinf$ be the theory containing $\To$ which expresses that the
value group is elementary equivalent to the additive group $\ZZ$ and
that the residue field is a pseudo-finite field. Each model of
$\Tinf$ is elementary equivalent with $k\llp t \rrp$ for some
pseudo-finite field $k$, because of either Ax, Kochen \cite{AK1},
\cite{AK3}, Er{\v s}ov \cite{Ersov}, Cohen \cite{Cohen}, Pas
\cite{Pas}, or others.

For any $d\in\NN_0$, let $\ZZ^{(d)}$ be the additive group
$\ZZ[r^{-1}]_{r\in I_d}$ with $I_d=\{r\in\NN_0\mid{\rm
gcd}(r,d)=1\}$ and $\NN_0=\{z\in\ZZ\mid z>0\}$. Let $\Tinfd$ be the
theory containing $\To$ expressing that the value group is
elementary equivalent with $\ZZ^{(d)}$ and that the residue field is
a pseudo-finite field.

If $k$ is a pseudo-finite field of characteristic zero, we denote by
$k\llp t \rrp^{(d)}$ the field
$$
\bigcup_{r\in I_d} k\llp t^{1/r}\rrp,
$$
where we take the union in a compatible way, that is, such that
$(t^{1/r})^m=(t^{1/r'})^{m'}$ for any integers $m,m'\geq0$,
$r,r'>0$, whenever $m'/r'=m/r$. Note that $k\llp t \rrp^{(d)}$ is a
model of $\Tinfd$ and that each model of $\Tinfd$ is elementary
equivalent with $k\llp t \rrp^{(d)}$ for some pseudo-finite field
$k$, because of similar Ax - Kochen - Er{\v s}ov principles as cited
above.

Let $\cT$ be as before.

\begin{definition}\label{drel}
A collection $P$ of relations $P_K$ on $K^n$ for each model $K$ of
$\cT$ which is uniformly definable by a formula in $\LPas$ is called
a \emph{definable relation over $\cT$}. Here, $n\in\NN_0$.
\end{definition}

\begin{definition}\label{dequi}
A collection $\sim$ of equivalence relations $\sim_K$ on $K^n$ for
each model $K$ of $\cT$ which is uniformly definable by a formula
$\psi(x,y)$ in $\LPas$, where $x$ and $y$ run over $K^n$, is called
a \emph{ definable equivalence relation over $\cT$}. Here,
$n\in\NN_0$.

\end{definition}

\begin{definition}\label{dfinite}
Let $\sim $ be a definable equivalence relation over $\cT$. If there
exists $N\in\NN_0$ such that, for any model $K$ of $\cT$, the number
of equivalence classes of $\sim_K$ is $\leq N$, then $\sim$ is
called \emph{finite over $\cT$}.

If $K$ is a model of $\cT$, we say that $\sim $ is \emph{finite over
$K$} if the number of equivalence classes of $\sim_K$ is finite.
\end{definition}

\begin{definition}\label{dtame}
A formula in $\Ltame$ with no quantifiers running over the valued
field sort is called a \emph{tame formula}.
\end{definition}

\begin{definition}\label{dtrel}
A definable relation $P$ over $\cT$ is called \emph{tame over $\cT$}
(resp.~\emph{tame over $K$}, for any given model $K$ of $\cT$), if
there exists a tame formula $\varphi(x)$ such that
$$
\cT\vdash \varphi(x)\ \leftrightarrow\ P(x)
$$
(resp.~$K\models \varphi(x)\ \leftrightarrow\ P(x)$).
\end{definition}

\begin{definition}\label{dtameim}
Let $\sim$ be a definable equivalence relation over $\cT$. Say
that \emph{the imaginaries of $\sim$ are tame over $\cT$}, if
there exists a tame formula $\psi(x,\xi,m)$, with $\xi$ a tuple of
residue field variables and $m$ a tuple of variables running over
$(G/nG)_{n=2,3,\ldots}$, such that
\begin{equation}\label{imtam}
\cT\vdash (\forall x)(\exists \xi)(\exists m)(\forall y)\big(x\sim
y \leftrightarrow \psi(y,\xi,m) \big).
\end{equation}
For  $K$ a model of $\cT$, say that \emph{the imaginaries of $\sim$
are tame over $K$} if there exists a $\psi$ as above (that may
depend on $K$) that satisfies the same condition but with
$\cT\vdash$ replaced by $K\models$.
\end{definition}

\begin{lem}\label{lfin}
Let $\sim$ be a definable equivalence relation over $\Tinf$. Suppose
that $\sim $ is finite over $k\llp t \rrp$  for each pseudo-finite
field $k$ of characteristic zero. Then $\sim$ is finite over
$\Tinf$.
\end{lem}
\begin{proof}
This follows by compactness in a standard way, cf.~the proof of
Lemma \ref{ltinf} for a less standard compactness argument.
\end{proof}

\begin{prop}\label{QEtame}
The theory $\Tinfd$ has elimination of valued field quantifiers
in $\Ltame$. Moreover, any $\LPas$-formula without free value group
variables is equivalent over $\Tinfd$ to a $\Ltame$-formula without
valued field quantifiers.
\end{prop}
\begin{proof}

We first prove that the collection of groups
$G=\ZZ^{(d)},G/2G,G/3G,\ldots$ has elimination of $G$-quantifiers in
the multisorted language $\LL_0$ consisting of $\LL_{\ord}$ for $G$,
together with the language of groups for each of the groups $G/nG$,
and the natural projection maps $\pi_{mn}:G/mG \to G/nG$ and
$\pi_n:G\to G/nG$ for $n$ dividing $m>0$.

It is enough to eliminate the quantifier $(\exists x)$ from a
formula $\varphi(y)$ of the form:
\begin{equation}\label{exists1}
(\exists x)\Big(\, \wedgem_i f_i(y)=K_ix \wedgem_j g_j(y)\not=L_jx
\wedgem_\ell \big( h_{\ell}(y)\sq_{\ell} M_\ell x\big)\wedgem
H_1(y)\wedgem H_2(\pi_m(x,y))\Big)
\end{equation}
with $y=(y_1,\ldots,y_n)$ running over $G^n$, $f_i,g_j,h_{\ell}$
linear homogeneous forms in $y$ over $\ZZ$, $m>0$ and the $K_i$,
$L_j$, $M_\ell$ integers, the $\sq_{\ell}$ either $<$ or no
condition, the $H_i$ formulas without $G$-quantifiers, and
$\pi_m(x,y)=(\pi_m(x),\pi_m(y_1),\ldots,\pi_m(y_n))$.
 By changing $H_1$ if necessary and because the order on each coset of $mG$ in $G$ is dense,
we may suppose that all the $\sq_{\ell}$ are no condition.
 Suppose first that at least one of the $K_i$ is nonzero. Then, again
changing $H_1$ if necessary, we may suppose that there are no
polynomials $g_j$. But then (\ref{exists1}) is equivalent with a
certain $G$-quantifier free formula $H_3(y)$.
 Suppose finally that there are no polynomials $f_i$. Then we may
suppose that there are no polynomials $g_j$ since the fibers of the
maps $\pi_m$ are infinite. Then clearly (\ref{exists1}) is
equivalent with a certain $G$-quantifier free formula $H_4(y)$. This
proves the $G$-quantifier elimination in $\LL_0$.

\par
Now we prove the statement for the theory $\Tinfd$. The theory
$\Tinfd$ has elimination of valued field quantifiers and of
$G$-quantifiers in the language $\LPas\cup\LL_0$ by \cite{Pas} and
by the above quantifier elimination result for $G$-quantifiers in
$\LL_0$. We have to show that we can remove the sort for $G$.
 Let $\varphi(x,\xi,\alpha)$ be a $\LPas\cup\LL_0$-formula without
valued field quantifiers and $G$-quantifiers, where $x$ are valued
field variables, $\xi$ are residue field variables, and $\alpha$
runs over $(G/nG)_{n=2,3,\ldots}$, but possibly containing the
symbol $\ord$. For $f_i$ polynomials over $\ZZ$, the condition
$$
\ord\, f_i(x)< \ord\, f_j(x)
$$
with possibly $f_j(x)=0$ and $f_i(x)\not=0$, is equivalent to the
condition
$$
\ac (f_i(x)+f_j(x) )=\ac (f_i(x)+2f_j(x) )=\ac (f_i(x)+3f_j(x))
$$
and one can rewrite conditions $ \ord\, f_i(x)\leq \ord\, f_j(x) $
and $ \ord\, f_i(x) = \ord\, f_j(x) $ similarly. This easily shows
that $\varphi(x,\xi,\alpha)$ is equivalent to a $\Ltame$-formula
without valued field quantifiers.
\end{proof}

\begin{lem}\label{ltinf}
Let $P$ be a definable relation over $\Tinf$. Suppose that $P$ is
tame over $k\llp t \rrp$  for each pseudo-finite field $k$ of
characteristic zero. Then $P$ is tame over $\Tinf$.
\end{lem}
\begin{proof}Let $x$ be the tuple consisting of the variables that occur freely in $P$.
For any tame formula $\psi(x)$, let $C_\psi$ be the sentence
$$
(\forall x) \left( P(x)\leftrightarrow \psi(x)\right).
$$
By the supposition of the lemma and by compactness,
$$
\Tinf\vdash C_{\psi_1}\vee\ldots\vee C_{\psi_n}
$$
for some formulas $\psi_1,\ldots,\psi_n$ and some $n$.
 For each $j=1,\ldots,n$, let $D_j$ be the sentence
$$
C_{\psi_j}\wedge \left(\neg C_{\psi_1}\wedge\ldots\wedge \neg
C_{\psi_{j-1}}  \right),
$$
where $\neg$ is the negation.
Because $D_j$ has no free variables, we see by elimination of valued
field quantifiers \cite{Pas} that $D_j$ is equivalent over $\Tinf$
with a sentence in the residue field language, and hence equivalent
over $\Tinf$ with a tame formula.

  Now let $\psi(x)$ be the formula
$$
\left(D_1\rightarrow \psi_1(x)\right)\wedge\ldots\wedge
\left(D_n\rightarrow \psi_n(x)\right).
$$
Then, $\psi$ is
equivalent over $\Tinf$
 with a tame formula and
$$
\Tinf\vdash \psi(x)\ \leftrightarrow\ P(x)
$$
by the construction of the proof, and hence, $P$ is tame over
$\Tinf$.
\end{proof}

\begin{lem}\label{ltinfd}
Let $\sim$ be a definable equivalence relation over $\Tinfd$.
Suppose that the imaginaries of $\sim$ are tame over $k\llp t
\rrp^{(d)}$ for each pseudo-finite field $k$ of characteristic zero.
Then the imaginaries of $\sim$ are tame over $\Tinfd$.
\end{lem}

\begin{proof}
Although the proof is similar to the proof of Lemma \ref{ltinf},
we give the details.

Suppose that $\sim$ is an equivalence relation in $n$ variables and
let $x=(x_1,\ldots,x_n)$ run over the valued field. For any tame
formula $\psi(x,\xi,m)$, with $\xi$ a tuple of residue field
variables and $m$ a tuple of variables running over
$(G/nG)_{n=2,3,\ldots}$, let $C_\psi$ be the sentence
$$
(\forall x)(\exists \xi)(\exists m)(\forall y)\left(x\sim y
\leftrightarrow \psi(y,\xi,m) \right).
$$
By the supposition of the lemma and by compactness,
$$
\Tinfd\vdash C_{\psi_1}\vee\ldots\vee C_{\psi_n}
$$
for some formulas $\psi_1,\ldots,\psi_n$ and some $n$. By taking
the tuples $\xi$ and $m$ big enough, we may suppose that the
$\psi_j$ have free variables included in $x,\xi,m$. For each
$j=1,\ldots,n$, let $D_j$ be the sentence
$$
C_{\psi_j}\wedge \left(\neg C_{\psi_1}\wedge\ldots\wedge \neg
C_{\psi_{j-1}}  \right),
$$
where $\neg$ is the negation.  Now let $\psi(x,\xi,m)$ be the
formula
$$
\left(D_1\rightarrow \psi_1(x,\xi,m)\right)\wedge\ldots\wedge
\left(D_n\rightarrow \psi_n(x,\xi,m)\right).
$$
Then, $\psi$ is
$\Tinfd$-equivalent with a tame formula, because of the same reason
as in the proof of Lemma \ref{ltinf},
and
$$
\Tinfd\vdash (\forall x)(\exists \xi)(\exists m)(\forall
y)\left(x\sim y \leftrightarrow \psi(y,\xi,m) \right)
$$
by the construction of the proof. Hence, the imaginaries of $\sim$
are tame over $\Tinfd$.

\end{proof}

\begin{prop}[Criterion]\label{propcrit}
Let $P$ be a definable relation over $\Tinf$ that can be defined by
$\theta(x)$ in $\LPas$ that is existential w.r.t.~the valued field
variables and the value group variables. Suppose that for each
pseudo-finite field $k$ of characteristic zero there exists
$d\in\NN_0$ such that for all $r\in\NN_0$ with ${\rm gcd}(r,d)=1$ we
have for all tuples $x$ over $k\llp t \rrp$
$$
k\llp t \rrp\models \theta(x) \quad \leftrightarrow\quad k\llp
t^{1/r}\rrp\models \theta(x).
$$
Then $P$ is tame over $\Tinf$.
\end{prop}
\begin{proof}
By Lemma \ref{ltinf} it suffices to prove that $P$ is tame over
$k\llp t \rrp$ for $k$ an arbitrary pseudo-finite field $k$ of
characteristic zero. Note that $P$ induces a relation on $k\llp t
\rrp^{(d)}$ which is defined by $\theta(x)$. By the hypothesis of
the Proposition,
$$
k\llp t \rrp\models \theta(x) \quad \leftrightarrow\quad k\llp t\rrp
^{(d)}\models \theta(x).
$$
Now apply Proposition \ref{QEtame} to $\theta$ to obtain a
$\Ltame$-formula $\psi$ without valued field and
$(G/nG)_{n=2,3,\ldots}$-quantifiers such that
$$
k\llp t\rrp ^{(d)}\models \theta(x) \quad \leftrightarrow\quad k\llp
t\rrp ^{(d)}\models \psi(x).
$$
We may suppose that only maps $\ord_n$ with $n$ dividing a power
of $d$ occur in $\psi$. For $n$ dividing a power of $d$,
the natural map from $G/nG$ for $k\llp t \rrp$ to $G/nG$ for $k\llp
t \rrp^{(d)}$ is an isomorphism.
Hence,
$$
k\llp t\rrp \models \psi(x) \quad \leftrightarrow\quad k\llp t\rrp
^{(d)}\models \psi(x),
$$
and thus, $P$ is tame over $k\llp t \rrp$.
\end{proof}

\begin{conjecture}\label{conj1}
Let $\sim$ be a definable equivalence relation over $\Tinf$. If
$\sim$ is finite and tame over $\Tinf$, then the imaginaries of
$\sim$ are tame over $\Tinf$.
\end{conjecture}

\begin{conjecture}\label{conj2}
Let $\sim$ be a definable equivalence relation over $\Tinfd$. If
$\sim$ is finite over $\Tinfd$, then the imaginaries of $\sim$ are
tame over $\Tinfd$.
\end{conjecture}

\begin{prop}\label{con1con2}
Conjecture \ref{conj2} implies Conjecture \ref{conj1}.
\end{prop}
\begin{proof}
This follows from the Transfer Lemma \ref{ltransfer}.
\end{proof}

\begin{remark}
If Conjecture \ref{conj2} is true, then Theorem \ref{mt1} follows
from it, from Lemma A, from our Criterion for Tameness
\ref{propcrit}, and from Proposition \ref{con1con2}, by a standard
ultraproduct argument (this way, one thus avoids the use of Lemmas B
and C of section \ref{seccoh}). These conjectures \ref{conj1} and
\ref{conj2} are related to elimination of imaginaries as studied in
model theory. In fact, Hrushovski \cite{Hnotes} recently proved
Conjecture \ref{conj2} by using techniques for obtaining elimination
of imaginaries, cf~\cite{UdiEI}.
\end{remark}

\begin{lem}[Transfer Lemma]\label{ltransfer}
Let $\sim$ be a definable equivalence relation over $\Tinf$, which
is tame over $\Tinf$, say, defined by a tame formula $\theta(x,y)$.
Then there exists a $d_0\in\NN_0$ such that for all multiples $d$ of
$d_0$ the formula $\theta(x,y)$ defines an equivalence relation
$\sim_d$ over $\Tinfd$ such that the following assertions hold.
 \item[(i)] For each model of the form $k\llp t \rrp^{(d)}$ of
 $\Tinfd$,
 the relation $\sim_d$ is the union of the relations $\sim$ on the
 subfields $k\llp t^{1/r}\rrp$.
 \item[(ii)] If $\sim$ is finite over $\Tinf$, then $\sim_{d}$ is finite over
 $\Tinfd$.
 \item[(iii)] If the imaginaries of $\sim_{d_0}$ are tame over $\cT_{\infty}^{(d_0)}$
then the imaginaries of $\sim_{d}$ are tame over $\Tinfd$ and the
imaginaries of $\sim$ are tame over $\Tinf$.
\end{lem}

\begin{proof}
Let $d_0$ be the product of all the $n$ such that the map $\ord_n$
occurs in $\theta$. Then (i) follows immediately. Suppose now that
$\sim$ is finite over $\Tinf$. Then, by definition, there exists
$N>0$ such that in any model of $\Tinf$ there are at most $N$
equivalence classes. Let $\theta_{N}(a_1,\ldots,a_{N+1})$ be the
tame formula $\vee_{i\not= j}\, \theta(a_i,a_j)$; if it holds for
all tuples $a_i$ then there are at most $N$ equivalence classes.
Since $\theta_N$ is a tame formula only involving $\ord_n$ with $n$
dividing $d_0$, $\theta_N(a_1,\ldots,a_{N+1})$ holds for any tuples
$a_i$ in the models $k\llp t \rrp^{(d)}$ and thus (ii) follows. Now
suppose that the imaginaries of $\sim_{d_0}$ are tame over
$\cT_{\infty}^{(d_0)}$, say, the
 tame
 formula $\psi(x,\xi,m)$ satisfies
\begin{equation}\label{tam1}
\cT_{\infty}^{(d_0)} \vdash (\forall x)(\exists \xi)(\exists
m)(\forall y)\big(x\sim_{d_0} y \leftrightarrow \psi(y,\xi,m)
\big).
\end{equation}
We may suppose that the only maps $\ord_n$ that occur in $\psi$ are
such that $n$ divides a power of $d_0$. Take a tuple $x$ over some
model $k\llp t \rrp^{(d)}$ of $\cT_{\infty}^{(d)}$. Since $k\llp t
\rrp^{(d)}\subset k\llp t \rrp^{(d_0)}$ and by (\ref{tam1}), we can
take $\xi,m$ such that
$$
k\llp t \rrp^{(d_0)} \models (\forall y)\big(x\sim_{d_0} y
\leftrightarrow \psi(y,\xi,m) \big).
$$
For $n$ dividing a power of $d_0$,
the natural maps from $G/nG$ for $k\llp t \rrp$ to $G/nG$ for $k\llp
t \rrp^{(d)}$ and from both these to $G/nG$ for $k\llp t
\rrp^{(d_0)}$ are isomorphisms.
Hence,
$$
k\llp t \rrp^{(d)} \models (\forall y)\big(x\sim_d y
\leftrightarrow \psi(y,\xi,m) \big)
$$
and if moreover $x\in k\llp t \rrp$ then also
$$
k\llp t \rrp \models (\forall y)\big(x\sim y \leftrightarrow
\psi(y,\xi,m) \big).
$$
In other words,
$$
k\llp t \rrp^{(d)} \models (\forall x)(\exists \xi)(\exists
m)(\forall y)\big(x\sim_d y \leftrightarrow \psi(y,\xi,m) \big)
$$
and
$$
k\llp t \rrp \models (\forall x)(\exists \xi)(\exists m)(\forall
y)\big(x\sim y \leftrightarrow \psi(y,\xi,m) \big),
$$
which finishes the proof of (iii) by Proposition \ref{ltinfd}.
\end{proof}

\section{Cohomological Lemmas}\label{seccoh}

When $K$ is a field, we denote by $K^a$ an algebraic closure of
$K$. For a field $K$ and a discrete topological group $G$ with a
continuous action of the Galois group $\Gal (K^a / K)$, we denote
the first cohomology set of $G$ by $H^1(K,G)$.

Moreover, if $L$ is a Galois extension of $K$, we will also
consider $H^1(L,G)$ and $H^1(L/K,G(L))$, where the last one is the
cohomology set of the group $G(L)$, of all elements in $G$ which
are fixed under $\Gal (K^a/L)$, with respect to the obvious action
of $\Gal (L/K)$ on $G(L)$.

For all these notions we refer to \cite[I \S 5 and III]{Serre1}.
Note that these cohomology sets are pointed sets : they are
equiped with a distinguished element, so that we can speak about
exact sequences.

\begin{prop}\label{prop1} Let $k$ be a pseudo finite field of characteristic zero, and
$G$ a linear algebraic group over $K := k((t))$. Then we have

a) There are only a finite number of fields between $K$ and $K^a$ of
any given finite degree over $K$.

b) $H^1(K,G)$ is finite.

c) If $G$ is semi-simple and simply connected, then $H^1(K,G) = 0$.
\end{prop}

\begin{proof}

a) Let $L$ be such a field of degree $n$ over $K$ and with
ramification index $e$. Then $L = k^\prime((t))(\sqrt[e]{at})$
with $[k^\prime:k] = n/e$ and $a \in k^\prime$. Thus $L \subset
k^{\prime \prime}((t)) (\sqrt[e]{t})$, where $k^{\prime \prime}$
contains all extensions of $k^\prime$ inside $k^a$ with degree
$\leq e$.

b) This follows from a) and a theorem of Borel and Serre, see
\cite[III.4.3 Th\'eor\`eme 4]{Serre1}.

c) A pseudo finite field has cohomological dimension $\leq 1$, see
e.g. \cite[Corollary 10.19]{Fried-Jarden}. But a theorem of
Bruhat-Tits \cite{Bruhat-Tits} asserts that $H^1(K,G) = 0$, when $G$
is semi-simple and simply connected over $K$, whenever $K$ is a
complete field with respect to a discrete valuation, whose residue
field has cohomological dimension $\leq 1$. (In the $p$-adic case
this is known as Kneser's Theorem.)

\bigskip
\noindent {\bf Lemma A.}\textit{ Let $k$ be a pseudo finite field
of characteristic zero, and $G$ a linear algebraic group over $K =
k((t))$. Then there exists an integer $d \geq 1$ having the
following property : If $L$ is any field extension of $K$ with $n
:= [L:K] < \infty$, and $\gcd(d,n) = 1$, then the restriction map
$$H^1(K,G) \stackrel{\rm res}{\longrightarrow} H^1(L,G)$$
is injective.}

\bigskip
\noindent {\it Proof.} By twisting (see \cite[I \S 5.3 and III \S
1.3]{Serre1}) with cocycles in a set of representatives of
$H^1(K,G)$, which is finite by Proposition \ref{prop1}, it
suffices to prove that the kernel of res is zero. There are the
following cases:

\bigskip
\noindent \textit{ \textbf{1) $G$ is connected and reductive.}}

\noindent Let $G^{ss}$ be the (semi-simple) derived group of $G$,
and $G^{sc}$ the universal covering group of $G^{ss}$. Let $T$ be
the connected component of the center of $G$. It is a torus. The
kernel $\Delta$ of the epimorphism $T \times G^{sc} \rightarrow G$,
induced by the multiplication map $T \times G^{ss} \rightarrow G$,
is finite and contained in the center of $T \times G^{sc}$ (see
\cite[p. 325]{Borel} and \cite[2.2.2.(3), p. 37]{Springer}).
 This yields an exact sequence
$$0 \rightarrow \Delta \rightarrow T \times G^{sc} \rightarrow G
\rightarrow 0.$$ By Proposition \ref{prop1} c), this induces a
commutative diagram of exact sequences
\begin{equation*}\xymatrix{
H^1(K,\Delta)  \ar[d]_{\rm res} \ar[r] & H^1(K,T) \ar[d]_{\rm res}
\ar[r] & H^1(K,G) \ar[d]_{\rm res} \ar[r] & H^2(K,\Delta)
\ar[d]_{\rm res}
\\
H^1(L,\Delta)   \ar[r] & H^1(L,T) \ar[r] & H^1(L,G)  \ar[r] &
H^2(L,\Delta).}
\end{equation*}

The cohomology sets of $T$ and $\Delta$ are abelian groups and for
these one can consider the corestriction maps. Let $d$ be the
least common multiple of the orders of the groups $\Delta$ and
$H^1(K,T)$. Let $L$ be a finite extension of $K$ whose degree $n$
is relatively prime to $d$. Then the first, second, and fourth
vertical arrows in the above diagram are injective, because
composing them with the corestriction morphisms yields
multiplication by $n$, which is bijective on elements that are
anihilated by $d$. Each element of $H^1(K,T)$ which is mapped to
zero in $H^1(L,G)$ belongs to the image of $H^1(K,\Delta)$,
because corestriction commutes with the most left horizontal
arrows in the above diagram. Straightforward diagram chasing now
shows that the third vertical arrow has indeed a trivial kernel.

\bigskip
\noindent \textit{\textbf{ 2) $G$ is connected.}}

\noindent Let $G_u$ be the unipotent radical of $G$. Because
$H^1(K,-)$ is zero on any twist of $G_u$, by \cite[III.2.1 Prop.
6]{Serre1}, we see that $H^1(K,G)$ injects into $H^1(K,G/G_u)$.
Apply now the previous case to the  connected reductive group
$G/G_u$ to obtain the desired result.

\bigskip
\noindent \textit{\textbf{ 3) $G$ is finite.}}

\noindent Since $H^1(K,G)$ is finite (by Proposition \ref{prop1}),
there exists a finite Galois extension $K^\prime$ of $K$ such that
the inflation map (which is always injective)
$$H^1(K^\prime / K,G(K')) \rightarrow H^1(K,G)$$
is bijective, and such that $G(K^\prime) = G$. Put $d =
[K^\prime:K]$. Let $L$ any finite extension of $K$ whose degree $n$
is relatively prime with $d$. Consider the commutative diagram

\begin{equation*}\xymatrix{
H^1(K'/K,G(K'))  \ar[d]_{\alpha} \ar[r]^{\quad\ \cong} & H^1(K,G)
\ar[d]_{\rm res}
\\
H^1(K'.L/L,G(K'.L))   \ar[r] & H^1(L,G)}
\end{equation*}

where $\alpha$ is the obvious natural map. Note that $\alpha$ is a
bijection because the natural map from $\Gal (K^\prime \cdot L/L)$
to $\Gal (K^\prime / K)$ is an isomorphism (since $K^\prime \cap L
= K$) and $G(K^\prime \cdot L) = G(K^\prime) = G$. Moreover, the
bottom horizontal map in the above diagram is injective, because
it is an inflation map. Thus res is indeed injective.

\bigskip
\noindent  \textit{\textbf{ 4) $G$ is any linear algebraic group.}}

\noindent Let $G_0$ be the identity component of $G$, and $E =
G/G_0$ the finite quotient. Let $K'$
 be a finite Galois extension of $K$ such that $E(K^\prime) =
E$. Let $d$ be a positive integer such that Lemma A is true for $G$
replaced by $G_0$ (Case 2), and for $G$ replaced by $E$ (Case 3),
and such that $d$ is divisible by $[K^\prime : K]$. For $L$ any
finite extension of $K$ whose degree $n$ is relatively prime with
$d$, one has $E(K) = E(L)$. Indeed, the natural map from $\Gal
(K^\prime \cdot L/L)$ to $\Gal (K^\prime/K)$ is an isomorphism since
$K^\prime \cap L = K$. The desired result is now obtained by
straightforward diagram chasing in the following diagram with exact
rows:

\begin{equation*}\xymatrix{
E(K)  \ar@{=}[d] \ar[r] & H^1(K,G_0) \ar[d]_{\rm res} \ar[r] &
H^1(K,G) \ar[d]_{\rm res} \ar[r] & H^1(K,E) \ar[d]_{\rm res}
\\
E(L)   \ar[r] & H^1(L,G_0) \ar[r] & H^1(L,G)  \ar[r] & H^1(L,E).}
\end{equation*}

This finishes the proof of Lemma A.
\end{proof}

\begin{remark}
From the above proof, and the material in the proof of Proposition
\ref{prop1}, it is clear that Lemma A remains true when $K$ is any
complete characteristic zero field with respect to a discrete
valuation, whose residue field has cohomological dimension $\leq 1$,
such that there are only a finite number of fields between $K$ and
$K^a$ of any given finite degree over $K$. When $K$ is a $p$-adic
field and $G$ connected, this was proved by a different method in a
paper of Sansuc \cite[Remarque 4.8.1]{Sansuc}.
\end{remark}

\bigskip
\noindent {\bf Lemma B.}\textit{ Let $K$ be a field of
characteristic zero, and $k$ a subfield of $K$. Let $G$ be a linear
algebraic group over $K$. Suppose that $G$ is obtained from a linear
algebraic group $G_k$ over $k$ by base change and denote by $G(k^a)$
the group of rational points on $G_k$ over the algebraic closure
$k^a$ of $k$ in $K^a$. Then the natural map
$$H^1 (K,G(k^a)) \rightarrow H^1(K,G)$$
is surjective.}

\begin{proof} We consider the following cases:

\bigskip
\noindent \textit{ \textbf{ Case 1:}} $G$ has a normal algebraic
subgroup $T$ defined over $k$ which is a torus, such that $E := G/T$
is finite.

\bigskip
\noindent Consider the following commutative diagram with exact
rows:

\begin{equation*}\xymatrix{
H^1(K,T(k^a))  \ar[d] \ar[r] & H^1(K,G(k^a)) \ar[d]_{\alpha} \ar[r]
& H^1(K,E) \ar@{=}[d]
\\
H^1(K,T)   \ar[r] & H^1(K,G) \ar[r] & H^1(K,E). }
\end{equation*}

\noindent Let $b \in H^1(K,G)$. We have to prove that $b$ is in
the image of $\alpha$. Let $c$ be the image of $b$ in $H^1(K,E)$,
and let $\bar c$ be a cocycle representing $c$. We denote by
$\Delta(\bar c)$ the image of $\bar c$ in $H^2(K,\, _{\bar
c}T(k^a))$, cf.~\cite[I \S 5.6]{Serre1}, where $_{\bar c}T(k^a)$
is obtained from $T(k^a)$ by twisting with the cocycle $\bar c$.

\bigskip
\noindent{\bf Claim 1.} The natural map from $H^2(K,\, _{\bar
c}T(k^a))$ to $H^2(K,\, _{\bar c}T)$ is injective. (Here $_{\bar
c}T$ is obtained from $T$ by twisting with $\bar c$).

\bigskip
\noindent We prove this claim later, and first proceed with the
proof of Case 1. The natural image of $\bar c$ in $H^2(K,\, _{\bar
c}T)$ is zero, because $c$ is the image of $b \in H^1(K,G)$.
Hence, Claim 1 implies that $\Delta (\bar c) = 0$ in $H^2(K,\,
_{\bar c}T(k^a))$. Thus $c$ belongs to the image of
$H^1(K,G(k^a))$ in $H^1(K,E)$, and there exists an $a \in
H^1(K,G(k^a))$ such that $\alpha(a)$ and $b$ have the same image
$c$ in $H^1(K,E)$.

\bigskip
Let $\bar a$ be a cocycle representing $a$. We can twist the above
diagram by the cocycle $\bar a$. Let $_{\bar a}G(k^a),\, _{\bar
a}G,\, _{\bar a} T(k^a),\, _{\bar a}T$ be obtained from $G(k^a),
G, T(k^a),T$, by twisting with $\bar a$. By \cite[I \S 5.5
Corollary 2]{Serre1} we see that $b$ is in the image of $H^1(K,
_{\bar a}T)$, under the natural map from $H^1(K,\, _{\bar a}T)$ to
$H^1(K,\, _{\bar a}G)$ composed with the canonical bijection
between $H^1(K,\, _{\bar a}G)$ and $H^1(K,G)$. Hence, to prove
that $b$ is in the image of $\alpha$, it suffices to show the
following claim:

\bigskip
\noindent {\bf Claim 2.} The natural map $H^1(K,\, _{\bar
a}T(k^a)) \rightarrow H^1(K,\, _{\bar a} T)$ is surjective.

\bigskip
\noindent We now prove Claim 2. Let $w$ be any element of the
abelian torsion group $H^1(K,\, _{\bar a}T)$ and let $n$ be the
order of $w$. Let $T_n$ be the kernel of the $n$-th power map on
$_{\bar a}T(k^a)$. We have the following commutative diagram with
exact rows:

\begin{equation*}\xymatrix{
H^1(K,T_n)  \ar@{=}[d] \ar[r] & H^1(K,\, _{\bar a} T(k^a))
\ar[d]_{\alpha} \ar[r]^{(\cdot)^n} & H^1(K,\, _{\bar a} T(k^a))
\ar[d]
\\
H^1(K,T_n)   \ar[r] & H^1(K,\, _{\bar a} T) \ar[r]^{(\cdot)^n} &
H^1(K,\, _{\bar a} T). }
\end{equation*}

Notice that $w$ belongs to the image of the first arrow in the
second row, because $w$ belongs to the kernel of the second arrow in
that row. This implies Claim 2. Next we turn to the proof of Claim
1.

\bigskip
\noindent Put $W =\, _{\bar c}T/\, _{\bar c}T(k^a)$. We have an
exact sequence
$$H^1(K,W) \rightarrow H^2(K,\, _{\bar c}T(k^a)) \rightarrow
H^2(K,\, _{\bar c}T).$$ Note that $W$ is a uniquely divisible
abelian group. Hence $H^1(K,W) = 0$. This finishes the proof of
Claim 1, and thus also of Case 1.

\bigskip
\noindent
 \textit{ \textbf{ Case 2:}}
 $G$ is obtained by base change to $K$ from a linear
algebraic group over $k$ whose identity component is reductive.

\bigskip
\noindent Let $T$ be a maximal torus of $G$ defined over $k$. Thus
$T$ is a Cartan subgroup of $G$ and $T$ has finite index in its
normalizer $N$. A result of Springer (cf.~\cite[III.4.3, Lemma
6]{Serre1}) asserts that the natural map
$$H^1(K,N) \rightarrow H^1(K,G)$$
is surjective. Applying Case 1 to $N$ finishes the proof of Case 2.

\bigskip
\noindent
 \textit{ \textbf{ Case 3:}}
 $G$ is obtained by base change to $K$ from any linear
algebraic group over $k$.

\bigskip
\noindent Let $G_u$ be the unipotent radical of the identity
component of $G$. Note that the natural map from $H^1(K,G)$ to
$H^1(K,G/G_u)$ is injective (cf.~the argument in the proof of Case
2 of Lemma A). Moreover the natural map
$$H^1(K,G(k^a)/G_u(k^a)) \rightarrow H^1(K,G/G_u)$$
is surjective by Case 2. Hence to prove Case 3, it suffices to
prove the following claim.

\bigskip
\noindent {\bf Claim 3.} The natural map
$$H^1(K,G(k^a)) \rightarrow H^1(K,G(k^a)/G_u(k^a))$$
is surjective.

\bigskip
\noindent We now turn to the proof of Claim 3. This is well known in
the special case that $k=K$, cf.~\cite[Lemma 1.13]{Sansuc}; for the
general case we need a different argument. We may suppose that $G_u
\ne \{1 \}$. Let $C$ be the center of $G_u$, then $C$ is unipotent
and $\dim C \geq 1$ (see e.g. \cite[ \S 17]{Humphreys}). Moreover,
$C$ is connected, because any unipotent linear algebraic group over
a field of characteristic zero is connected, see chapter 3, section
2, Corollary 2 of Theorem 1 of \cite{OnishV}.

\noindent By induction on $\dim G$, it suffices to prove that the
natural map
$$H^1(K,G(k^a)) \rightarrow H^1(K,G(k^a)/C(k^a))$$
is surjective.

\noindent Using Proposition 41 in Chapter I.5.6 of \cite{Serre1}, we
see that it suffices to prove that $H^2(K,\, _{\bar e}C(k^a)) = 0$
for each cocycle $\bar e$ with values in $G(k^a)/C(k^a)$, where
$_{\bar e}C$ is obtained from $C$ by twisting with $\bar e$. But
this is clear because the abelian group $C(k^a)$ is uniquely
divisible, since it is unipotent (thus admitting a composition
series with successive quotients isomorphic to $k^a,+$). This
terminates the proof of Lemma B.
\end{proof}

\bigskip
\noindent {\bf Lemma C.}\textit{ Let $k$ be a field of
characteristic zero, and $G$ a linear algebraic group over $K =
k((t))$. Suppose that $G$ is obtained from a linear algebraic group
$G_k$ over $k$ by base change. Let $k^a$ be the algebraic closure of
$k$ in $K^a$ and let ${\cal O}$ be the valuation ring of $K^a$,
i.e., the integral closure of $k[[t]]$ in $K^a$. Denote $G_k(k^a)$
by $G(k^a)$ and $G_k({\cal O})$ by $G({\cal O})$. Then the map
$$H^1(K,G({\cal O})) \rightarrow H^1 (K,G(k^a)),$$
induced by reduction modulo the maximal ideal of ${\cal O}$, is
injective.}

\begin{proof}
 Let $I$ be the kernel of the reduction map
$G({\cal O}) \rightarrow G(k^a)$. We have an exact sequence
$$H^1(K,I) \rightarrow H^1(K,G({\cal O})) \rightarrow
H^1(K,G(k^a)).$$ Hence, it suffices to prove that $H^1(K,\, _cI) =
0$ for each 1-cocycle $c$ with values in $G({\cal O})$. Here $_cI$
denotes the twist of $I$ by $c$.

\bigskip
\noindent Let $L$ be any finite Galois extension of $K$, over
which $c$ is defined. We are going to prove that $H^1(L/K,\, _c
I(L)) = 0$, where $I(L) = I \cap G(L)$. This implies that
$H^1(K,\, _cI) = 0$.

\bigskip
\noindent We look at the filtration $$I(L) = I_1 \triangleright
I_2 \triangleright I_3 \triangleright \cdots \triangleright I_j
\triangleright \cdots,$$ with $$I_j := \{ g \in G({\cal O}) \cap
G(L) \mid g \equiv 1 \bmod \pi^j \},$$ for $j = 1,2,\cdots,$ where
$\pi$ is a generator for the maximal ideal of ${\cal O} \cap L$.

\bigskip
\noindent Note that $I_j/I_{j+1}$ is an abelian group, for each $j
\geq 1$. This is easily verified by identifying $G$ with a subgroup
of $GL_n$ for a suitable $n$. Moreover, the abelian group
$I_j/I_{j+1}$ is uniquely divisible. Indeed this follows from
Hensel's Lemma, because the map $G \rightarrow G : g \mapsto g^m$
induces multiplication by $m$ on the tangent space of $G$ at $1$.
Thus $H^1(L/K,\, _c I_j/\, _c I_{j+1}) = 0$, for all $j \geq 1$.

\bigskip
\noindent The exact sequence
$$0 \rightarrow \frac{_c I_j}{ \, _cI_{j+1}} \rightarrow \frac{_cI_1}{_cI_{j+1}} \rightarrow
\frac{_cI_1}{_cI_j} \rightarrow 0$$

induces an exact sequence
$$H^1(L/K,\, _cI_j/\, _c I_{j+1}) \rightarrow H^1(L/K,\, _cI_1/\, _cI_{j+1})
\rightarrow H^1(L/K,\, _cI_1/\, _cI_j).$$ Using induction on $j$,
we conclude that
\begin{equation}\label{star}
H^1(L/K,\, _cI_1/\, _c I_j) = 0 \mbox{ for all } j \geq 1.
\end{equation}
Next we turn to the proof that $H^1(L/K,\, _cI(L)) = 0$. Let $a =
(a_\sigma)_{\sigma \in \Gal (L/K)}$ be any 1-cocycle of $\Gal
(L/K)$ with values in $_cI(L) =\, _cI_1$. From (\ref{star}) it
follows that for each $j \geq 1$ there exists $b_j \in\, _cI_1$
such that for all $\sigma \in \Gal (L/K)$ we have
\begin{equation}\label{starstar}
a_\sigma \equiv b^{-1}_j \sigma(b_j) \quad \bmod \quad  \pi^j.
\end{equation}
 Let
$\kappa$ be the ramification index of $L$ over $K$, and assume
$\kappa \mid j$. We can choose a basis for ${\cal O \cap L}$ over
$k[[t]]$, and write the  components of the tuple $b_j$ in terms of
that basis (considering the coefficients as unknowns), and obtain
in this way a system of polynomial equations whose solvability in
$k[[t]]/t^{j/\kappa}$ is equivalent with the existence of $b_j$
satisfying (\ref{starstar}).
 Applying Greenberg's Theorem \cite{Greenberg1} to this
system of equations, we see that there exists $b \in\, _cI_1$ such
that
$$a_\sigma = b^{-1} \sigma(b)$$
for all $\sigma \in \Gal  (L/K)$. This finishes the proof that
$H^1(L/K,\, _cI(L)) = 0$. Lemma $C$ is now proven.
\end{proof}

\section{Definability of cohomology}\label{secdef}

For $\cL$ a language and $a=(a_i)_i$ a tuple (or a set) of elements
of an $\cL$-structure, $\cL(a)$ denotes the language $\cL$ with
extra constant symbols for the $a_i$.

\begin{lem}\label{finited}
Let $k$ be a pseudo-finite field of characteristic zero. Let $d$ be
an integer. Then the field $K:=k\llp t \rrp^{(d)}$ has only finitely
many field extensions of any given finite degree inside an algebraic
closure $K^a$ of $K$.
\end{lem}
\begin{proof}
The lemma follows from the fact that a field extension $L$ of degree
$n$ of $K$ is contained in $K_{m,n}:=(k_m\llp t
\rrp^{(d)})[t^{1/n}]$, with $k_m$ the unique field extension of $k$
of degree $m$ and with $m$ big enough so that $k_m$ contains $k_n$
and all $n$-th roots of $1$. We leave the details to the reader.
\end{proof}

Let $K$ be $k\llp t \rrp^{(d)}$.  Let $e$ be a positive integer.
Choose $n$ such that all $d^e$-th roots of $1$ are in $k_n$, with
$k_n$ the unique degree $n$ extension of $k$. Let $M$ be the
finite Galois extension of $K$
$$M:=k_n\llp t \rrp^{(d)}.$$

Let $L$ be the finite Galois extension of $K$
\begin{equation}\label{L}
L:=M[t^{1/d^e}].
\end{equation}

Note that any finite field extension of $K$ is contained in a
field $L$ as in (\ref{L}).
\begin{lem}\label{defGal}
Let $a_i\in k$, $i=0,\ldots,n-1$, be the coefficients of an
irreducible polynomial $f_a(x):=x^n+\sum_{i=0}^{n-1} a_ix^i$ over
$k$. Then, the field $k_n$, the group ${\rm Gal}(L/K)$, and its
action on $k_n$ are $\Ltame(a)$-definable, that is, they are
isomorphic to a $\Ltame(a)$-definable field and group which acts
definably on the field.

For $i=0,\ldots,n-1$, let $a'_i\in \cO_K$ be such that $a_i'$ lies
above $a_i$ and let $a_n'\in \cO_K$ be such that $\ac(a_n')=1$ and
such that the image of $\ord(a_n')$ in $G/d^eG$ is a generator of
$G/d^eG$. Write $a=(a_1\ldots,a_{n-1})$ and
$a'=(a'_1,\ldots,a'_{n-1},a'_{n})$. Then, the field $L$ with the
action of ${\rm Gal}(L/K)$ is $\Ltame(a',a)$-definable, uniformly in
$a'$, that is, it is given by formulas in which the tuple $a'$ may
occur but which is independent of the choice of $a'$.
\end{lem}
\begin{proof}
If $\xi$ is a zero of $f_a$, $k_n$ is the vector space $k^n$ with
multiplicative structure induced by the isomorphism $k^n\to
\bigoplus_{i=0}^{n-1}\xi^ik$, which is independent of the choice
of $\xi$. Hence, $k_n$ is $\Ltame(a)$-definable.

The Galois group of $k_n$ over $k$ is cyclic of order $n$, say,
with generator $\sigma$, and each of its powers $\sigma^\ell$
corresponds to a matrix $B_\ell=(b_{ij})$ in $GL_n(k)$ by
\begin{equation}
\sigma^\ell(\xi^j)=\sum_{i=0}^{n-1}b_{ij}\xi^i.
\end{equation}
Since $\Gal(k_n/k)$ is commutative, each matrix $B_\ell$ is
independent of the choice of the zero $\xi$. Moreover, the matrix
subgroup $B:=\{B_\ell\}_\ell$ of $GL_n(k)$ is independent of the
choice of generator $\sigma$.

Clearly, the field $M$ with the action of ${\rm Gal}(M/K)$ is
$\Ltame(a,a'_1,\ldots,a_{n-1}')$-definable, uniformly in
$a'_1,\ldots,a_{n-1}'$.

Now we use $a'=(a'_1,\ldots,a'_{n})$ to define $L$. For $a'_n$ with
$\ac(a_n')=1$ and such that the image of $\ord(a_n')$ in $G/d^eG$ is
a generator of $G/d^eG$, there exists by Hensel's Lemma $\xi\in K$
and an integer $b$ with $\gcd(b,d)=1$ such that
$$
t^b=\xi^{d^e}a'_n.
$$
Hence, $M[t^{1/d^e}]$ is the same field as $M[(a_n')^{1/d^e}]$
within a fixed algebraic closure of $M$. Thus the field $L$ is
$\Ltame(a,a')$-definable, uniformly in $a'$.

Let $\mu_{d^e}\in k_n$ be a primitive $d^e$-th root of unity. The
Galois group of $L$ over $M$ is cyclic of order $d^e$ and isomorphic
to the multiplicative group $(\mu_{d^e})$ generated by $\mu_{d^e}$.
 Since ${\rm Gal}(L/K)$ is a semidirect product $(\mu_n) \rtimes B$,
it is clearly $\Ltame(a)$-definable. The action on $k_n$ is clearly
also $\Ltame(a)$-definable. Thus the field $L$ with the action of
${\rm Gal}(L/K)$ is clearly $\Ltame(a,a')$-definable, uniformly in
$a'$.
\end{proof}

For $x_0$ in $X(K)$ (or in $X(\cO_K)$), let $H_{x_0}$ be the
stabilizer of $x_0$, considered as a linear group over $K$, resp.~as
a group scheme over $\cO_K$. Following Serre \cite[III.4.4]{Serre1},
there is an injection ${\rm Cor}_{x_0}$ from the orbits in $X(K)$
under the action of $G(K)$ into $ H^1(K,H_{x_0}(K^a))$. Namely,
${\rm Cor}_{x_0}(\mbox{orbit of }x)$ with $x \in G(K)$ is defined to
be the class of the cocycle
\begin{equation}\label{orbittau}
\sigma\mapsto \tau^{-1}\sigma(\tau),
\end{equation}
where $\tau\in G(K^a)$ is such that $x=\tau(x_0)$.

Fix $x_0\in X(\cO_K)$ and let $\bar H_{x_0}$ be the reduction of
$H_{x_0}$ modulo the maximal ideal of $\cO_K$. Since
$H^1(K,H_{x_0}(K^a))$ is finite by \cite[III.4.3 Th\'eor\`eme
4]{Serre1} and Lemma \ref{finited}, we may assume that $L$ is big
enough (and still of the form (\ref{L})) so that
$H^1(K,H_{x_0}(K^a))$ is naturally isomorphic to
$H^1(L/K,H_{x_0}(L))$. Hence we can consider the diagram
\begin{equation}\label{diagr}\xymatrix{
X(K) \ar[r]^{p} & X(K)/G(K) \ar[d]_{{\rm Cor}_{x_0}}\\
H^1(L/K,H_{x_0}(\cO_L))  \ar[d]_{i_{x_0}} \ar[r]^{\pi_{x_0}}  &
H^1(L/K,H_{x_0}(L))
\\
H^1(L/K,\bar H_{x_0}(k_n)),}
\end{equation}
with $\cO_L$ the valuation ring of $L$, and $p$, $\pi_{x_0}$, and
$i_{x_0}$ the natural maps.

From now on until the end of the proof of Theorem \ref{mt4} we
suppose that $F=\QQ$. The treatment for general number fields $F$ is
completely similar but requires coefficients from $F$ in all the
valued field and residue field languages and respective definitions
considered in section \ref{sec:setting}; we leave it to the reader
to carry this out (cf.~the comments preceding Theorem \ref{mt5}
or~\cite{DL}, \cite{CLoes} for adding coefficients to a language).

To say that a certain subset of $X(K)$ is $\Ltame(c)$-definable we
work with a finite cover with affine charts of $X$, defined over
$\QQ$,
 cf.~\cite{DL} or \cite{CLoes}.

\begin{prop}\label{def4}
Let $x$ be in $X(K)$. Suppose that there exists $x_0\in X(\cO_K)$
such that $\pi_{x_0}$ is surjective and such that $i_{x_0}$ is
injective. Then there exists $c\in k^m$ for some $m$ such that the
orbit of $x$ under the action of $G(K)$ is $\Ltame(c)$-definable.
(The point is that no direct reference to the point $x$ can be made
in $\Ltame(c)$.)
\end{prop}
\begin{proof}
Let $E(x_0)$ be the condition on $x_0$ that $x_0\in X(\cO_K)$,
that $\pi_{x_0}$ is surjective, and that $i_{x_0}$ is injective.
 Let $x_0$ satisfy $E$.
 By Serre \cite[III.4.4]{Serre1} we know that ${\rm Cor}_{x_0}$ is injective.
 It is enough to prove that the condition
 \begin{equation}\label{ec1}
(\exists x_0)\, \left(E(x_0) \wedge i_{x_0}\, \pi^{-1}_{x_0}\,
{\rm Cor}_{x_0}\, p\, (x) = i_{x_0}\, \pi^{-1}_{x_0}\, {\rm
Cor}_{x_0}\, p\, (y) \right)
 \end{equation}
on $y\in G(K)$ is $\Ltame(c)$-definable for some $c\in k^m$, since
it cuts out the orbit of $x$. Take $a\in k^n$ as in Lemma
\ref{defGal}. Let $D(a')$ be the condition on $a'=(a'_i)_{i=0}^{n}$
that $a_i'\in\cO_K$ lies above $a_i$ for $i<n$, and that $a_n'\in
\cO_K$ is such that $\ac(a_n')=1$ and such that the image of
$\ord(a_n')$ in $G/d^eG$ is a generator of $G/d^eG$. By Lemma
\ref{defGal}, the field $L$ with the action of ${\rm Gal}(L/K)$ on
it is $\Ltame(a,a')$-definable, uniformly in $a'$. Also $k_n$,
$\cO_L$, the projection $\cO_L\to k_n$, and the action of ${\rm
Gal}(L/K)$ on $k_n$ and on $\cO_L$ are $\Ltame(a,a')$-definable,
uniformly in $a'$. Hence, $H_{x_0}(L)$, $\bar H_{x_0}(k_n)$,
$H_{x_0}(\cO_L)$, the natural maps $H_{x_0}(\cO_L)\to \bar
H_{x_0}(k_n)$ and $H_{x_0}(\cO_L)\to H_{x_0}(L)$, and the action of
${\rm Gal}(L/K)$ on $H_{x_0}(L)$, $\bar H_{x_0}(k_n)$, and on
$H_{x_0}(\cO_L)$ are $\Ltame(a,a',x_0)$-definable, uniformly in $a'$
and $x_0$. Since ${\rm Gal}(L/K)$ is finite, that two given
cocycles, $a_\sigma$, $b_\sigma$ representing elements of one of the
occurring $H^1(L/K,\cdot)$ are cohomologous is
$\Ltame(a,a',x_0)$-definable over the graphs of $a_\sigma$ and
$b_\sigma$ (more precisely, over all the entries of all the tuples
in these graphs), uniformly in $a'$, $x_0$, and in the graphs of
$a_\sigma$ and $b_\sigma$. One then readily checks that the
condition $E(x_0)$ is $\Ltame(a,a')$-definable, uniformly in $a'$,
that is, it is given in the charts by formulas in which $a'$ may
occur but which is independent of the choice of $a'$.

Let $a_\sigma$ be a cocycle representing an element of
$H^1(L/K,\bar H_{x_0}(k_n))$ which lies in $i_{x_0}\,
\pi^{-1}_{x_0}\, {\rm Cor}_{x_0}\, p\, (x)$, for some $x_0$
satisfying $E$. By identifying cocycles with their graph we may
use them in formulas. In particular, (the graph of) $a_\sigma$ is
$\Ltame(c)$-definable for some $c\in k^m$. Let $F(a'_\sigma,x_0)$
be the property on $a'_\sigma$ that $a'_\sigma$ is a cocycle
representing an element of $H^1(L/K, H_{x_0}(\cO_L))$ that is
mapped to the class of $a_\sigma$ under $i_{x_0}$. Likewise, let
$G(a'_\sigma,x_0,y)$ be the property on $a'_\sigma,y$ that the
class of $a'_\sigma$ is mapped to ${\rm Cor}_{x_0}p y$ under
$\pi_{x_0}$.
 As before, the conditions $F$ and $G$ are
$\Ltame(a,c,a',x_0)$-definable, uniformly in $a'$ and $x_0$.

It is clear that
\begin{equation}\label{ec2}
(\exists a')(\exists x_0)(\exists a'_\sigma)\, \left(D(a')\wedge
E(x_0)\wedge F(a'_\sigma,x_0)\wedge G(a'_\sigma,x_0,y)\right)
\end{equation}
is equivalent to (\ref{ec1}). Moreover, (\ref{ec2}) is
$\Ltame(a,c)$-definable, where $a,c$ are the tuples in $k$ as
constructed above. This proves the Proposition.

\end{proof}

\section{Proof of the main results}\label{secproof}
We prove the following slight generalization of Theorem \ref{mt1}.
\begin{theorem}\label{mt4}Let $U\subset X$ be an affine open, defined over $F$.
Then there are finitely many regular functions $f_i:U\to \AA^1_F$
and integers $N>0$ and $d>0$ such that, for any $K\in \cC_N$ and any
$x\in X(K)$, the set $G(K)(x)\cap U(K)$ depends only on
$\ac(f_i(x))$ and a statement in the language of groups on the
$\ord(f_i(x))\bmod d$, that is, interpreted in the value group
modulo $d$, and where $\ord(0)\equiv 0\bmod d$ by convention.
\end{theorem}

\begin{proof}
Recall that we suppose till the end of the proof that  $F=\QQ$. The
proof for general $F$ is completely similar but requires
coefficients from $F$ in all the valued field and residue field
languages and respective definitions considered in section
\ref{sec:setting}; we leave it to the reader to carry this out
(cf.~the comments preceding Theorem \ref{mt5} or~\cite{DL},
\cite{CLoes} for adding coefficients to a language).

For any field $K$ over $\QQ$ let $\sim_K$ be the equivalence
relation on $X(K)$ such that two points are equivalent for
$\sim_K$ if and only if they lie in the same orbit under the
action of $G(K)$. Let $\sim$ be the equivalence relation over
$\To$ whose interpretation is $\sim_K$ for any $K$ which is a
model of $\To$.

Let $k$ be a pseudo-finite field of characteristic zero.
 Since $X$ is absolutely irreducible, there exists $x_0\in
X(k)$.\footnote{This is the only place where the absolutely
irreducibility of $X$ is used; see Remark \ref{remirr}.}
 Lemma A implies the conditions of Criterion \ref{propcrit} for
$\sim_{k\llp t\rrp}$ and for some $d$. This can be seen by applying
Lemma $A$ to the group $H_{x_0}$ with $x_0\in X(k)$ and by using the
injectivity and naturality of ${\rm Cor}_{x_0}$, as defined in
(\ref{orbittau}) following~\cite[III.4.4]{Serre1}. Hence, the
conclusion of Criterion \ref{propcrit} holds, that is, $\sim$ is
tame over $\Tinf$.

Now let $K$ be $k\llp t\rrp^{(d)}$ with $d>0$. A straightforward
adaptation of the proof of \cite[I.2.2, Prop.~8]{Serre1} yields that
$$
H^1(K,A)=\lim\limits_{\longrightarrow} H^1(k\llp
t\rrp[t^{1/r}],A),
$$
when $A$ is a discrete topological group with a continuous action of
the Galois group $\Gal (K^a/k\llp t\rrp )$, and where $r$ runs over
all integers with $\gcd(r,d)=1$, directed by divisibility, since the
absolute Galois group of $K$ is isomorphic to the projective limit
of the abolute Galois groups of the $k\llp t\rrp[t^{1/r}]$. (Note
that formation of equivalence classes commutes with taking direct
limits.)

Since direct limits preserve surjectivity and injectivity (in the
category of sets), Lemma's B and C imply that, for $x_0\in
X(k)\subset X(\cO_K)$ and for well chosen $L$, the map $\pi_{x_0}$
is surjective and that the map $i_{x_0}$ is injective in diagram
(\ref{diagr}). This holds for all $d>0$.

 Take $x\in U(K)$.
 Now apply Proposition \ref{def4}. One finds some $c\in k^m$ for
some $m$ such that the orbit of $x$ in $X(K)$ under the action of
$G(K)$ is $\Ltame(c)$-definable. By the quantifier elimination
result \ref{QEtame}, this orbit is given by a $\Ltame(c)$-formula
$\psi(c,y)$ where the only quantifiers that may occur in $\psi$ run
over the residue field. Since $\sim_K$ is finite by \cite[III.4.3
Th\'eor\`eme 4]{Serre1} and Lemma \ref{finited}, one obtains that
the imaginaries of $\sim_K$ are tame over $K$ by combining the
formulas $\psi(c,y)$ for all orbits. Namely, let $A_i$, with
$i=1,\ldots,\ell$ for some $\ell$, be the equivalence classes of
$\sim_K$, associate to each $A_i$ a formula $\psi_i(c_i,y)$, for
some $c_i$ running over $k^m$ for some $m$ as above, take a tuple of
variables $\xi=(\xi_0,\xi_1,\ldots,\xi_m)$ running over $k^{1+m}$,
 and let $\psi$ be the formula
$$
\vee_{i=1}^\ell \left(\psi_i(\xi_1,\ldots,\xi_m,y)\wedge
\xi_0=i\right).
$$
Then
\begin{equation}\label{imtam2}
K\models  (\forall x)(\exists \xi)(\exists m)(\forall y)\big(x\sim_K
y \leftrightarrow \psi(y,\xi,m) \big).
\end{equation}
In other words, the imaginaries of $\sim_K$ are tame over $K$.

By Lemma \ref{ltinfd}, the imaginaries of $\sim$ are tame over
$\Tinfd$ and this holds for any $d>0$.

By taking an appropriate $d$ and since $\sim$ is tame over $\Tinf$,
the Transfer Lemma \ref{ltransfer} asserts that the imaginaries of
$\sim$ are tame over $\Tinf$. Now Theorem \ref{mt4} and thus Theorem
\ref{mt1} follow for some $d$, by an ultraproduct argument.
\end{proof}

Write $\Ltame(F)$ for the language $\Ltame$ together with
coeficients from $F$ in the valued field and residue field sort.
We will prove Theorem \ref{mt5} which is a slight generalization
of Theorem \ref{mt3}. For the notion of definable subassignments
and definable morphisms we refer to \cite{DL} and \cite{CLoes}. By
$\Ltame(F)$-definable subassignments or $\Ltame(F)$-definable
morphisms we mean definable subassignments, resp.~ definable
morphisms, which are also definable in the language $\Ltame(F)$;
this means that they are given by finitely many
$\Ltame(F)$-formulas in affine charts defined over $F$. In
particular, one can take $K$-rational points on any
$\Ltame(F)$-definable subassignment for any field $K$ over $F$
which carries a $\Ltame(F)$-structure and similarly for
$\Ltame(F)$-definable morphisms.

\begin{theorem}\label{mt5} Write $X'$ for $X\otimes_F \Spec (F\llp t\rrp)$. Let $f:X'\to \AA^1_{F\llp t\rrp}$
and $g:X'\to \AA^1_{F\llp t\rrp}$ be $\Ltame(F)$-definable
morphisms. Let $\omega$ be a volume form on an affine open in $X$.
Let $W$ be a $\Ltame(F)$-definable subassignment
of the functor which sends a field $E$ over $F$ to $X'(E\llp
t\rrp)$.
For each $K\in \cC_1$ and for each $x$ in $X(K)$, under the
condition of integrability for all $s>0$, consider the orbital
integral
\begin{equation}\label{int5}
I_{K,x}(s):=\int_{G(K)(x) \cap W(K)} |g_K|\, |f_K|^s\, |\omega|_K,
 \end{equation}
with $|\omega|_K$ the measure on $X(K)$ associated to $\omega$,
and $g_K$, $f_K$, and $W_K$ the $K$-interpretations of $f$ and
$W$. If for some $x$ and $K$ this is not integrable, put $
I_{K,x}(s):=0$ for this $K$ and $x$.

Then, the conclusion of Theorem \ref{mt3} and its addendum holds.
\end{theorem}

\begin{proof}
We give two proofs. (In fact, the proof holds for tame integrals,
that is, integrals over a tame domain of a tame integrand,
cf.~\cite{Denefdegree} and \cite{Pas} for the notion of tame
integrals and Remark \ref{remt} for an extension of this notion;
that the domain of the integral (\ref{int5}) is tame follows from
Theorem \ref{mt4}.)

Firstly, one obtains Theorem \ref{mt5} from Theorem \ref{mt4} by
taking a suitable embedded resolution of singularities with normal
crossings and calculate the integral on the resolution space,
cf.~\cite{Denefdegree}.
 Secondly, one obtains Theorem \ref{mt5} from
Theorem \ref{mt4} by using the elementary method of Pas \cite{Pas}
to calculate $I_{K,x}(s)$ uniformly in $K\in\cC_N$ for $N$ big
enough, cf.~\cite{Pas}.
\end{proof}

\begin{remark}
Note that the rationality of $I_{K,x}(s)$ as in Theorem \ref{mt5}
and the fact that only finitely many poles with bounded multiplicity
can occur already follows from the rationality results for motivic
integrals in \cite{CLoes}, by enriching the Denef-Pas languages with
constant symbols for a point in each of the orbits for the theory
$\Tinf$; this approach is based on Denef-Pas cell
decomposition~\cite{Pas}.
\end{remark}

\begin{remark}\label{remirr}
The fact that $X$ is absolutely irreducible is not really needed for
Theorems \ref{mt1}, \ref{mt3}, its addendum, and the results of this
section. Indeed, instead of asking that $X$ is absolutely
irreducible, it is enough to ask that at least one irreducible
component over $F$ of $X$ is absolutely irreducible. This implies
that there exists a point in $X(k)$ for any pseudo-finite field $k$
over $F$ and this is all that the absolute irreducibility of $X$ was
used for.
\end{remark}

\begin{remark}\label{remt}
In fact, in Theorem \ref{mt5}, one can consider still more general
integrals for which the same conclusions hold. Indeed, the
integrand in (\ref{int5}) may be a finite sum of terms of the form
\begin{equation}\label{eqt}
\left(\prod_{j=1}^\ell \ord (h_{j,K})\right) |g_K|\, |f_K|^s\,
|\omega|_K,
\end{equation}
where the $|g_K|$, $|f_K|^s$, and $|\omega|_K$ are as in the
theorem, the $h_j:X'\to \GG_m$ are $\Ltame(F)$-definable morphisms
with $K$-interpretation $h_{j,K}$, and $\ord$ is the order on
$K^\times$ extended by zero. This is so because an integral with
integrand (\ref{eqt}) over a tame domain can be rewritten as a tame
integral over more variables with an integrand as in (\ref{int5}).
\end{remark}

\subsection*{Acknowledgment}
The authors are grateful to Mikhail Borovoi for several helpful
remarks and to the referee for careful remarks.\\

{\small During the realization of this project, the first author was
a postdoctoral fellow of the Fund for Scientific Research - Flanders
(Belgium) (F.W.O.) and was supported by The European Commission -
Marie Curie European Individual Fellowship with contract number HPMF
CT 2005-007121.}

\bibliographystyle{amsplain}
\bibliography{anbib}

\end{document}